
\magnification=\magstep1
\vsize=22truecm
\input amstex
\documentstyle{amsppt}
\leftheadtext{E. Makai, Jr., H. Martini}
\rightheadtext{Unique local determination of convex bodies}
\topmatter
\title Unique local determination of convex bodies
\endtitle
\author E. Makai, Jr., H. Martini\endauthor
\address MTA Alfr\'ed R\'enyi Institute of Mathematics,
\newline
H-1364 Budapest, Pf. 127, Hungary
\newline
{\rm{http://www.renyi.mta.hu/\~{}makai}}
\newline
\vskip.1cm
Fakult\"at f\"ur Mathematik, Technische Universit\"at Chemnitz,
\newline
D-09107 Chemnitz, Germany
\endaddress
\email makai.endre\@renyi.mta.hu, martini\@mathematik.tu-chemnitz.de\endemail
\thanks 
\endthanks
\keywords questions of Barker-Larman and Yaskin-Zhang,
unique determination of convex bodies, quermassintegrals, radial function,
star-shaped set, characterizations of symmetry\endkeywords
\subjclass {\it 2010 Mathematics Subject Classification.}
Primary: 52A2. Secondary: 52A38\endsubjclass
\abstract 
Barker and Larman asked the following. Let $K' \subset {\Bbb{R}}^d$ be a
convex
body, whose interior contains a given convex body $K \subset {\Bbb{R}}^d$, and
let, for all supporting hyperplanes $H$ of $K$, the 
$(d-1)$-volumes of the intersections $K' \cap H$ be given. 
Is $K'$ then uniquely determined? 
Yaskin and Zhang asked the analogous question when,
for all supporting hyperplanes $H$ of $K$,
the $d$-volumes of the ``caps'' cut off from $K'$ by $H$ are given.
We give local
positive answers to both of these questions, 
for small $C^2$-perturbations of $K$,
provided the boundary of $K$ is $C^2_+$. In both cases, $(d-1)$-volumes
or $d$-volumes can be replaced by $k$-dimensional
quermassintegrals for $1 \le k \le d-1$ or for $1 \le k \le d$, respectively.
Moreover, in the first case we can admit, rather than hyperplane sections, 
sections by $l$-dimensional affine planes, where $1 \le k \le l \le d-1$. 
In fact, here
not all $l$-dimensional affine subspaces are needed, but only a small
subset of them (actually, a $(d-1)$-manifold), for unique local
determination of $K'$.
\endabstract
\endtopmatter\document

\head \S 1 Introduction \endhead

Barker and Larman [BL], p. 81, Conjecture 2, posed the following question.
Let $K' \subset {\Bbb{R}}^d$ be a convex
body whose interior ${\text{int}}\,K'$
contains a given convex body $K \subset {\Bbb{R}}^d$, and
let, for all supporting hyperplanes $H$ of $K$, the areas, i.e.,
$(d-1)$-volumes of the intersections $K' \cap H$ be given. 
Is then $K'$ uniquely determined? The paper [BL]
investigated only the
case when $K$ is the unit ball $B^d$ and obtained several partial results
to this question, for which we refer to their paper or, for 
some of them, cf. below. 
(There arises the question what happens 
if we replace hyperplane sections by
sections by any affine $l$-planes $Y \subset {\Bbb{R}}^d$, for 
fixed $l \in [1,d-2]$, and we know the
$l$-volumes of the intersections $K' \cap Y$ for all $Y$ supporting
$K$  --- i.e., $Y$ intersects $K$ but not its interior. 
However, this would follow easily if we knew the result for
hyperplanes. In fact, knowledge of all these $l$-volumes 
would determine the
intersection of $K'$ with any translate $Z+x$ 
of any linear $(l+1)$-subspace $Z$ of ${\Bbb{R}}^d$, such that $(Z+x) \cap 
({\text{int}}\,K) \ne \emptyset $ --- hence also would determine $K'$.)

The recent paper of Ryabogin-Yaskin-Zvavitch [RYZ] 
repeated this question in p. 332 as Question 8,
and its special case where $K=B^d$ in p. 331 as Question 7. The same was done
in [YZh], in Problem 1.2 and Problem 1.3.
Question 19 in [RYZ], 
p. 335 (that is, as they observe, equivalent to their Question 20, p.
336) is a special case of their Question 7. Namely, it asks the following.
If $B^d \subset
{\text{int}}\,K'$ and the areas of the intersections of $K'$ with any two
parallel different tangent
hyperplanes of $B^d$ are equal, is $K'$ then $0$-symmetric?
[RYZ] and [YZh] also give good overviews
about results concerning questions of this
type, as well as several new questions. 

Yaskin-Zhang [YZh], Problem 1.4, posed the following question.
Let $K' \subset {\Bbb{R}}^d$ be a convex
body whose interior contains a given convex body $K \subset {\Bbb{R}}^d$, and
let, for all supporting hyperplanes $H$ of $K$, the 
$d$-volumes cut off from $K'$ by the supporting hyperplanes $H$ of $K$
be given. Is then $K'$ uniquely determined?

Barker-Larman [BL], Conjecture 1, which is repeated for the special case
$K=B^d$ in [RYZ], p. 335, Question 18, is
of similar type as those treated in this paper, but our methods do not yield
its local solution. This question is: let $K \subset {\text{int}}\,K'$, where 
$K,K' \subset {\Bbb{R}}^d$ are convex bodies with $K \subset {\text{int}}\,K'$, 
and let us have for each hyperplane $H$
supporting $K$ that $K' \cap H$ is centrally symmetric. Is then $K'$
centrally symmetric (for $K=B^d$ $0$-symmetric), or even an ellipsoid? 
Some considerations about this question see in \S 2, Problem. If $K$ is
replaced by a point, e.g., $K=\{ 0 \} $, and {\it{for some $l \in [2, d-1]$
and for each affine $l$-plane $P \ni 0$ we have that
$K' \cap P$ is centrally symmetric, then either $K$ is $0$-symmetric, or $K$ is
an ellipsoid.}} This is called the {\it{``false centre theorem''}}, 
cf. \cite{BL}, p. 80. (There only the case $d-l=1$ is mentioned, but then,
fixing $l$,
a trivial induction for $d$ proves the above mentioned general result. 
Some further references to this problem are
Aitchison-Petty-Rogers [APR], Larman [L], Montejano and Morales-Amaya [MM-A],
Larman and Morales-Amaya [LM-A], V. Soltan [So],
J. Jer\'onimo-Castro and T. B. McAllister [J-CM].)

We cite only some theorems. The first one is
due to Santal\'o [Sa], see also [BL],
Theorem 1: {\it{for $d=2$ and $K$ the unit disc with centre $0$, 
a concentric circle $K'$ is uniquely determined. I.e., 
if all above sets}} $K' \cap H$ (chords of $K'$ tangent to $K$)
{\it{have constant length, then $K'$ is a circle with centre $0$}}. In fact,
[Sa] proved this {\it{also for convex curves on}} $S^2$. 
This was reproved once more by
Gorkavyy-Kalinin [GK], who also {\it{proved the analogous statement on the
hyperbolic
plane}}, and gave some further planar situations in which their proof works.

The next theorem is due to V. Yaskin [Y], Theorem: {\it{for $K$ the unit ball
in ${\Bbb{R}}^d$, among convex polytopes, $K'$ is uniquely
determined.}} For dimension $2$, a more general result
was earlier proved by G. Xiong, Y.-w. Ma,
W.-s. Cheung [XMCh], for $K$ $0$-symmetric and ${\text{bd}}\,K$, the boundary
of $K$, ``nice''.
I.e., {\it{if for two convex polytopes (polygons) $K'_1,K'_2$ containing $K$
in their interiors all above sets 
$K'_1 \cap H$ and $K' _2 \cap H$}} (for $H$ being tangent hyperplanes of $K$)
{\it{have equal areas (lengths), then}} $K'_1=K'_2$.

Next we cite [BL], Theorem 4, which 
{\it{answers the above question in the positive sense for $K=B^d$, 
when instead of hyperplane sections
one considers sections with $l$-dimensional affine planes 
tangent to $B^d$ and their $l$-volumes, 
for fixed $l \in [1, d-2]$}}. However, observe 
that the supporting affine $l$-planes to any convex body $K \subset {\Bbb{R}}^d$
form an $\left( (l+1)(d-l)-1 \right)$-manifold 
(for $l=d-1$ a $(d-1)$-manifold), while 
$$
{\text{for }} 1 \le l \le d-2 {\text{ we have }} (l+1)(d-l)-1 > d-1 \,. 
\tag 0.1
$$

\newpage

The unknowns (values of the radial function for all
$u \in S^{d-1}$) form a $(d-1)$-manifold.
That is, intuitively one has for $l \in [1, d-2]$
``much more equations, namely  
$\left( (l+1)(d-l)-1 \right)$-manifold many'' (the 
values of the $l$-volumes of all
the sections of $K'$ by supporting affine $l$-planes of $K$) ``than unknowns,
namely $(d-1)$-manifold many'' (values of the radial function of $K'$). This 
``explains solvability'' of this problem for $K=B^d$ and $l \in [1, d-2]$,
while the original problem with $K=B^d$ and 
$l=d-1$ is unsolved. In the statement of
our Theorem 1 we eliminate this discrepancy between the dimensions of the
respective manifolds.

The following theorem is due to [F] and [Ga83]: {\it{if for two distinct
interior points $p_1,p_2$
of $K'$ and any hyperplane $H$ containing any
one of these points the area of the
intersection $K' \cap H$ is given, then $K'$ is uniquely determined}}. Observe
that here the hypothesis implies knowledge of the even parts 
$\left( \varrho _i
^{d-1}(u)+\varrho _i^{d-1}(-u) \right) /2$ of the functions $
\varrho _i^{d-1}(u)$, where $\varrho_i (u)$ is the radial function of
$K' - p_i$, cf. [Gr], Theorem 5.6.3. That is, intuitively, we have ``half
information'' both for $p_1$ and $p_2$, which together uniquely determine
$K'$. So here the heuristics works. As observed in [RYZ], p. 332, this
statement is a variant [BL], p. 81, Conjecture 1 (cf. the first paragraph of the
Introduction), when we replace the convex body $K$ by a non-degenerate
segment. What happens
if we replace $K$ by a non-empty compact convex set of fixed dimension in
$[2,d-1]$? In particular, what happens if $K$ is replaced by 
a ball of the respective dimension?

The questions of Barker-Larman and Yaskin-Zhang
seem to be difficult even for the plane with
$K=B^2$. 
Rotating the sections $K' \cap H$ by moving $u$ (the unit outer normal of $H$)
in $S^1$ in the positive sense, and differentiating with respect to $u$ 
(for $\partial K$ 
being $C^2_+$), we obtain some equations containing the values and
the first and second
derivative values of some function at two different points, and unicity of the
solution of the equation ought to be proved. So this way does not seem to
lead to a solution.

We remark that an analogous question was settled by [MM] locally about
the unit ball, under suitable smoothness hypotheses. 
That question was the following: {\it{if the $(d-2)$-dimensional
surface area --- or some other quermassintegral}} (or intrinsic volume)
{\it{of smaller positive dimension --- 
of the intersection of a convex body $K \subset {\Bbb{R}}^d$ with 
hyperplanes $u^{\bot } + \lambda u$, for any fixed $u \in S^{d-1}$, attains its
maximum, e.g., at $\lambda =0$, is $K$ then $0$-symmetric}}? 
Here $u^{\bot }$ {\it{is the linear $(d-1)$-subspace orthogonal to}} $u \in
S^{d-1}$. \cite{MM}, Theorem, {\it{settled this
question, under suitable smoothness hypotheses, 
``locally'', in the positive sense, close to the unit ball}}. More
generally, {\it{hyperplanes can be substituted by sections with 
$l$-dimensional affine
planes}} (and then we consider quermassintegrals of dimension in $[1,l]$),
{\it{and still the analogous statement holds}}.
{\it{The cases of \,$l$-volumes of sections with $l$-dimensional affine planes,
for $l=1$ and $2 \le l \le d-1$, have been
solved earlier in the positive sense}}, not only locally, but {\it{in
their 
original form}}, by Hammer [H], Theorem 1, and Makai-Martini-\'Odor [MM\'O],
Corollary 3.2, respectively. Even, similarly to our
Theorems 1 and 4, it suffices to consider some family ${\Cal{L}}_{l+1}$ of
linear $(l+1)$-subspaces whose union is ${\Bbb{R}}^d$, and assume the
maximality property only among translates of linear 

\newpage

$l$-subspaces lying in
some $L_{l+1} \in {\Cal{L}}_{l+1}$ (cf. the first paragraph of the proof of
our Theorem 1 and Remark 3). We will give in this paper,
in the same spirit, ``local solutions'' to the questions of
Barker-Larman and Yaskin-Zhang, cf. Theorems 1 and 4 below.

The questions of Barker-Larman and Yaskin-Zhang 
can be answered locally, close to any convex
body of class $C^2_+$, not only to $B^d$.
Let $K \subset {\Bbb{R}}^d$ be a $C^2_+$ 
convex body, with $0 \in {\text{int}}\,K$,
and let $K^t \supset K^0=K$ for $t \in [0,1]$ be a small $C^2$-perturbation
of $K$. Here we consider
$K^t$ to be given by its radial function
$\varrho ^t (u)$, for $u \in S^{d-1}$,
and we investigate $(\partial / \partial t) \varrho ^t (u)|_{t=0}$,
which is of course everywhere non-negative.
Suppose that, for all tangent hyperplanes of $K$, either 
the ``asymptotical behaviour'' of the areas of their
intersections with $K^t$, or the ``asymptotical behaviour'' 
of the volumes of the
``caps'' cut off from $K^t$ by them, for $t \to 0$, are given.
Then ``in first order'' the approximation $K^t$ is uniquely determined.

I.e., suppose that there is a $C^2$-deformation
$K^t$ of
$K$, with $K^t \supset K^0=K$, 
for parameter values $t \in [0,1]$ --- i.e., $[0,1] \times
S^{d-1} \ni (u,t) \mapsto \varrho ^t(u) \in (0, \infty )$ 
is a $C^2$ function. 
Suppose that we know either the ``asymptotical behaviour'' of the
areas ($(d-1)$-volumes) of the intersection 
of each tangent hyperplane of $K$ with $K^t$, or the ``asymptotical behaviour''
of the $d$-volume of the
``cap'' cut off from $K^t$ by each tangent hyperplane of $K$, 
for $t \to 0$. (The exact meaning of ``asymptotical behaviour'' will be given
in Theorems 1 and 4.) Then the
first partial derivative of the radial function of 
$K^t$ with respect to $t$ at $t=0$, for each $u \in S^{d-1}$, 
is uniquely determined, cf. our Theorems 1 and 4. 

In both cases we may replace $(d-1)$-volumes and $d$-volumes by
quermassintegrals (or intrinsic volumes)
of any lower positive dimension $k$, and the analogous statements
hold, cf. our Theorems 1 and 4. 

Like in [MM], Theorem, 
and [MM\'O], Corollary 3.2, in our Theorem 1 we may allow,
rather than sections by hyperplanes tangent to $K$, also sections 
with affine planes of lower, but positive dimension $l \,\,( \ge k$),
tangent to $K$, and the analogous statement
holds, cf. our Theorem 1. This is a local positive answer to Question 1 of
[YZh].
Here, however, recall from above that all tangent
affine $l$-planes to $K$,
which yield our equations about the ``asymptotical behaviours'',
form an $\left( (l+1)(d-l)-1 \right)$-manifold, 
while the unknowns (values of the partial derivative of the radial function
with respect to $t$, 
for $t=0$, and all $u \in S^{d-1}$) form a $(d-1)$-manifold,
so that we have ``much more equations than unknowns'' (cf. \thetag{0.1} above). 
The particular case of our 
Theorem 1, with ${\Cal{L}}_{l+1}$ being the family of
all linear $(l+1)$-subspaces
of ${\Bbb{R}}^d$ which contain some fixed linear $l$-subspace of
${\Bbb{R}}^d$, uses only
``$(d-1)$-manifold many'' tangent affine $l$-planes to $K$ (cf. the paragraph
before our Theorem 1), for the
``$(d-1)$-manifold many'' unknowns,
and still has a positive answer. Thus also here the heuristics works.

Concerning convex bodies, we will use the standard notations, cf. [Sch]. We
denote the
norm of a vector $x \in {\Bbb{R}}^d$ by $\| x \| $. 
A {\it{convex body in\,}} ${\Bbb{R}}^d$ is a compact convex set in ${\Bbb{R}}^d$
with non-empty interior. The {\it{boundary}} and 
{\it{interior of a set $X \subset {\Bbb{R}}^d$}} are denoted
by ${\text{bd}}\,X$ and ${\text{int}}\,X$, respectively.
$V(\cdot )$ will denote ($d$-dimensional) {\it{volume}}. 

\newpage

We write ${\text{conv}}\,X$ and ${\text{lin}}\,X$
for the {\it{convex hull}} and {\it{linear hull of a set}} 
$X \subset {\Bbb{R}}^d$, respectively. 
The {\it{unit ball}} and the {\it{unit 
sphere of\,}} ${\Bbb{R}}^d$ are denoted by $B^d$ and $S^{d-1}$, respectively.
For the {\it{volume of}} $B^d$ we write $\kappa _d$.
The {\it{quermassintegrals}} of non-empty compact convex sets $K$ (cf. [Sch],
Ch. 4) are denoted by $W_k(K)$, for $0 \le k \le d$.
Following [Sch], \S 4.2, we write $V_k(K) := \kappa _{d-k} ^{-1}
{d \choose k } W_{d-k}(K)$, which is called the {\it{$k$'th 
intrinsic volume of}} $K$, for $0 \le k \le d$.
Then, by [Sch], \S 4.2, for $r \in [0, \infty )$ we have
$V(K + r B^d ) = \sum _{i=0} ^d \kappa _{d-k} V_k(K) r^{d-k}$.
The intrinsic volumes $V_k(K)$ are monotonous, are positively homogeneous
of degree 
$k$, and are continuous in the Hausdorff metric. Moreover, they remain
unchanged if ${\Bbb{R}}^d$ is embedded in some higher dimensional Euclidean
space, and we consider $V_k(K)$ as the $k$'th intrinsic volume of $K$
considered as a subset of the higher dimensional Euclidean space ([Sch], \S
4.2, p. 210; in fact, they are characterized by this property among constant
multiples of $W_{d-k}$).
By a {\it{star-shaped set}}, or {\it{star-shaped hull of a set}} 
we mean a star-shaped set, or star-shaped hull of a set with respect to $0$.
The {\it{radial function}} of a compact
star-shaped set $X \subset {\Bbb{R}}^d$ is 
$\varrho ( \cdot ) : S^{d-1} \to [0, \infty )$, defined by $\varrho (u)
:= \max \{ r \in [0, \infty ) \mid ru \in X \} $.
A convex body $K \subset {\Bbb{R}}^d$ 
is $C^2_+$ if its boundary is a $C^2$ submanifold
of ${\Bbb{R}}^d$, with everywhere positive Gauss curvature.
We write {\it{area}} for $(d-1)$-volume. 
If a convex body $K \subset {\Bbb{R}}^d$ is smooth and strictly convex, 
and $Y$ is a tangent affine plane of $K$ of some dimension, 
then we will write $\{ y \} :=K \cap Y$, and 
the {\it{origin in}} $Y$ is chosen to be $y$.
The Dupin indicatrix of a $C^2_+$ convex body $K \subset 
{\Bbb{R}}^d$, at some $x \in {\text{bd}}\,K$, 
lying in  the tangent hyperplane $H$ of $K$ at $x$,
is obtained in the following way. Let in some fixed rectangular 
coordinate system with $x$ the origin and
$H$ the hyperplane given by $x_d = 0$, and with
$K$ lying above $H$, ${\text{bd}}\,K$ have a local representation
$x_d=f(x_1, \ldots , x_{d-1})$. (Then $f(0, \ldots , 0)=0$ and, for each
$i \in [1, d-1]$, also $f_{x_i}(0, \ldots , 0)=0$.) Then the 
{\it{Dupin indicatrix of $K$ at $x$}}, in the chosen rectangular coordinate
system, is $ \{ (x_1, \ldots , x_{d-1} ) \in H \mid
(1/2) \sum _{i,j = 1} ^{d-1} f_{x_i x_j}(0, \ldots , 0) x_i x_j = 1 \}
$.


\head \S 2 Theorems \endhead

Let $d \ge 2$, and $1 \le k \le l \le d-1$ be any integers. Let
$K \subset {\Bbb{R}}^d$ be a $C^2_+$ 
convex body, containing the origin in its
interior, with radial function 
$\varrho ( \cdot ) : S^{d-1} \to (0, \infty )$. Let, for
$t \in [0,1]$, the convex body $K^t$ be a one-parameter deformation of $K$, 
with radial function $\varrho ^t( \cdot ) :  S^{d-1} \to (0, \infty )$, with
$K^t \supset K^0=K$, 
and with $[0,1] \times S^{d-1} \ni (t,u) \mapsto \varrho ^t(u) \in (0, \infty
)$ being a $C^2$ function. 
Then $(\partial \varrho ^t(u) / \partial t)|_{t=0}$ is, for each $u \in
S^{d-1}$, non-negative.
Let ${\Cal{L}}_{l+1}$ denote a  
family of linear $(l+1)$-subspaces of ${\Bbb{R}}^d$,
whose union is ${\Bbb{R}}^d$. 

An example for ${\Cal{L}}_{l+1}$ is the family of all linear $(l+1)$-subspaces
of ${\Bbb{R}}^d$ which contain some fixed linear $l$-subspace of
${\Bbb{R}}^d$. Observe that this example forms a $(d-l-1)$-manifold. Then all
affine $l$-planes, tangent to $K$, of these linear $(l+1)$-subspaces form a
manifold of dimension $(d-l-1)+(l+1)-1=d-1$. 
For this example we have in the following Theorem 1
``$(d-1)$-manifold many equations'' (one for each $Y$ in Theorem 1, 
namely for $\lim _{\varepsilon \to 0}V_k(K^{\varepsilon }\cap L) / \varepsilon
^{k/2}$), for ``$(d-1)$-manifold 

\newpage

many
unknowns'' (values of $(\partial \varrho ^t(u) / \partial t)|_{t=0}$, for
all $u \in S^{d-1}$), 
which heuristically says that, according to
the dimension of the manifolds in question, we have ``as many
equations as unknowns''. 

\proclaim{Theorem 1}
Assume the hypotheses described before this theorem.
Then the $k$'th intrinsic volume (e.g., the $l$-volume)
of the intersection of $K^{\varepsilon }$ with each tangent affine 
$l$-plane $Y$ of $K$, 
divided by $\varepsilon ^{k/2}$, tends for $\varepsilon \to 0$
to a finite, non-negative limit, depending on $Y$. This limit as a
function of $Y$, taken only for those $Y$'s which
lie in some linear $(l+1)$-subspace of
${\Bbb{R}}^d$ (namely in ${\text{\rm{lin}}}\,Y$) 
belonging to ${\Cal{L}}_{l+1}$, 
uniquely determines $(\partial \varrho ^t(u) / \partial t)|_{t=0}$ for each $u
\in S^{d-1}$. 
\endproclaim

Thus in the hypothesis we have information about the
``asymptotical behaviour'' of the intrinsic volumes of the intersections, 
and we have a conclusion about the ``asymptotical 
behaviour'' of the radial functions, for $\varepsilon \to 0$. The same
holds for the intrinsic volumes of the ``caps'' cut off from $K'$ in
Theorem 4.

The statement for $1 \le 
l \le d-2$ is a relative of [BL], Theorem 4. There a global
uniqueness was proved, but only for $K=B^d$, and with 
all tangent affine $l$-planes of $K$. 

Also in [BL], Theorem 5, for $d$ odd,
only a countably
infinite set of values of $\varepsilon $ are used (sections by hyperplanes
intersecting $B^d$, with distances from $0$ in a countably infinite set).
In our Theorem 1 this would ``approximately'' 
correspond to a series of $\varepsilon $'s converging to $0$ (actually, to
different such sequences for different $Y$'s) --- however, this does not
substantially change
the statement of our Theorem 1. The same holds for our Theorem 4.


The following Corollary 2 is a local
version of [RYZ], p. 335, Question 19, in the same
sense, as Theorem 1 is the local version of the question of
Barker-Larman [BL], p. 81, Conjecture 2.

\proclaim{Corollary 2}
Let all the 
hypotheses of Theorem 1 hold. Additionally, let $T : {\Bbb{R}}^d \to 
{\Bbb{R}}^d$ be a linear isometry, such that $TK=K$. Let the
$k$'th intrinsic volume (e.g., the $l$-volume)
of the intersection of $K^{\varepsilon }$ with each
tangent $l$-plane $Y$ of $K$, 
lying in some linear $(l+1)$-subspace of \,${\Bbb{R}}^d$ 
(namely in ${\text{\rm{lin}}}\,Y$) 
belonging to ${\Cal{L}}_{l+1}$, 
and with the tangent $l$-plane $TY$ of $K$,
lying in some linear $(l+1)$-subspace of \,${\Bbb{R}}^d$ 
(namely in $T({\text{\rm{lin}}}\,Y)$)
belonging to $T({\Cal{L}}_{l+1})$,
divided by $\varepsilon ^{k/2}$, tend for $\varepsilon \to 0$ to the same
limit. Then $(\partial \varrho ^t(u) / \partial t)|_{t=0}$ 
has equal values for $u$ and $Tu$. In
particular, for $Tu=-u$ the function
$(\partial \varrho ^t(u) / \partial t)|_{t=0}$ is an even function of $u \in
S^{d-1}$.
\endproclaim

 
{\bf{Remark 3.}}
Suppose the hypotheses of Theorem 1.  
Let, for $Y$ a tangent affine $l$-plane of $K$, 
$E^Y$ denote the Dupin indicatrix of $K \cap {\text{lin}}\,Y$ at $y$ with $\{
y \} = K \cap Y$. For unicity of $E^Y$ we suppose that the lengths of the
vectors in the rectangular coordinate system in ${\text{lin}}\,Y$ are equal to
their lengths in the original polar coordinate system in ${\Bbb{R}}^d$.
Then {\it{in Theorem 1 and Corollary 2
we can replace $V_k( \cdot )$ by any functions}} 
$F^Y( \cdot )$ (there being
no compatibility conditions between the functions $F^Y( \cdot )$
for different $Y$'s) {\it{which have the following properties}}.
\newline
1) {\it{They are defined on $[0, \infty ) \times N^Y$, where $N^Y$ is 
some neighbourhood of the set}}\,\, ${\text{conv}}\,E^Y$ --- 
considered in the set of convex bodies in $Y$
with the topology of the 

\newpage

Hausdorff metric --- {\it{and have 
values in $[0, \infty )$, but are positive on $(0, \infty ) \times N^Y$.}}
\newline
2) {\it{They are positively homogeneous of degree $k$.}} 
\newline
3) {\it{They are monotonous.}}
\newline

In particular, $F^Y( \cdot )$ (even defined for all non-empty compact convex 
subsets of $Y$)
can be chosen as a mixed volume $V( (\cdot ),
\ldots , (\cdot ), K_{k+1}, \ldots , K_l)$, with $ (\cdot )$
occurring $k$ times, and $K_{k+1}, \ldots , K_l \subset Y$ being
any non-empty compact convex sets. (For mixed volumes cf. [Sch], Chapter 5.)
But there are many other possibilities,
e.g., the $k$'th intrinsic volume of the minimum volume circumscribed, or
maximum volume inscribed ellipsoid, etc.

\vskip.1cm


In the next theorem, a {\it{cap cut off from $K^t$ by a tangent hyperplane 
$H$ of
$K$}} is the intersection of $K^t$ and of the closed half-space $H^+$
of ${\Bbb{R}}^d$,
bounded by $H$ and not containing $K$. (Of course, here
it makes no sense to consider affine planes of lower dimension.)

\proclaim{Theorem 4}
Let all the hypotheses of Theorem 1 hold, with $l:=d-1$, but with $1 \le k \le
d$.
Then the $k$'th intrinsic volume (e.g., the $d$-volume)
of the cap cut off from $K^{\varepsilon }$ 
by each tangent hyperplane $H$ of $K$,
divided by $\varepsilon ^{k/2}$ for $1 \le k \le
d-1$, or by $\varepsilon ^{(d+1)/2}$ for $k=d$, respectively,
tends for $\varepsilon \to 0$ to a finite non-negative limit, depending on $H$.
This limit, as a function of $H$, uniquely determines
$(\partial \varrho ^t(u) / \partial t)|_{t=0}$. 
\endproclaim

Observe that for $d=1$ (with $K, K'$ then being segments with $K \subset
K'$) the function mapping $x \in {\text{bd}}\,K$ to the volume
(length) of the cap cut off from $K'$ by $\{ x \} $ 
trivially uniquely determines $K'$. So we have not only a local, but a global
solution. 
Therefore $d \ge 2$ is assumed in the hypotheses of Theorem 4.
The same holds also for the following Corollary 5.


\proclaim{Corollary 5}
Let all the 
hypotheses of Theorem 4 hold. Additionally, let $T : {\Bbb{R}}^d \to 
{\Bbb{R}}^d$ be a linear isometry such that $TK=K$.
Let the $k$'th intrinsic volume (e.g., the $d$-volume)
of the caps cut off from $K^{\varepsilon }$ by each
tangent hyperplane $H$ of $K$, and by the tangent hyperplane $TH$ of $K$,
divided by $\varepsilon ^{k/2}$ for $1 \le k \le d-1$, or by $\varepsilon
^{(d+1)/2}$ for $k = d$, respectively,
tend for $\varepsilon \to 0$ to the same limit.
Then $(\partial \varrho ^t(u) / \partial t)|_{t=0}$ has 
equal values for $u$ and $Tu$. In
particular, for $Tu=-u$, the function
$(\partial \varrho ^t(u) / \partial t)|_{t=0}$ is an even function of $u \in
S^{d-1}$.
\endproclaim


{\bf{Remark 6.}}
Suppose the hypotheses of Theorem 4, and let, for $H$ a tangent hyperplane of
$K$, $E^H$ be as in Remark 3 (there $E^Y$ is written).
Then {\it{in Theorem 4 and Corollary 5
we can replace $V_k( \cdot )$ by any functions}} 
$G^H( \cdot )$ (there being
no compatibility conditions between the functions $G^H( \cdot )$
for different $H$'s)
{\it{which have the following properties}}.
\newline
1) {\it{They are defined
on $[0, \infty ) \times {\tilde N}^H$, where ${\tilde N}^H$ is 
some neighbourhood of the set}}\,\, ${\text{conv}}E^H$ --- considered 
in the set of non-empty compact convex sets in ${\Bbb{R}}^d$,
with the topology of 
the Hausdorff metric --- {\it{and have values in $[0, \infty )$, but are
positive on $(0, \infty ) \times {\tilde N}^H$.}}
\newline
2) {\it{They are positively homogeneous of degree $k$.}}
\newline
3) {\it{They are monotonous.}}
\newline
4) {\it{They are continuous, from the restriction to $[0, \infty ) \times
{\tilde{N}}^H$ of the 
topology of the 

\newpage

Hausdorff metric of non-empty
compact convex sets in ${\Bbb{R}}^d$ to $[0,\infty )$.}}
\newline
5) {\it{If $k=d$, then still we have that $G^H$ is invariant under volume 
preserving affinities.}}

In particular, $G^H( \cdot )$ (even defined for all non-empty compact convex 
subsets of ${\Bbb{R}}^d$)
can be chosen as a mixed volume $V( (\cdot ),
\ldots , (\cdot ), K_{k+1}, \ldots , K_d)$, with $ (\cdot )$
occurring $k$ times, and $K_{k+1}, \ldots , K_d \subset {\Bbb{R}}^d$ being
any non-empty compact convex sets. 
But there are many other possibilities,
e.g., the $k$'th intrinsic volume of the minimum volume circumscribed, or
maximum volume inscribed ellipsoid, etc.


{\bf{Problem.}}
The question of Barker-Larman [BL], p. 81,
Conjecture 1, repeated in the special case
$K=B^d$ in [RYZ], p. 335, Question 18, is the following.
Let $K \subset {\text{int}}\,K'$, where 
$K,K' \subset {\Bbb{R}}^d$ are convex bodies,
and let us have for each hyperplane $H$ supporting
$K$ that $K' \cap H$ is centrally symmetric. Is then $K'$
centrally symmetric (for $K=B^d$ $0$-symmetric), or even an ellipsoid? 
(Of course, also here one could replace hyperplane sections by
sections
by any affine $l$-planes supporting $K$, for fixed $l \in [2,d-2]$, and suppose
their central symmetry. However, a
positive answer to \cite{BL}, p. 81,
Conjecture 1, would easily imply the analogous
statement for each $l \in [2,d-2]$. In fact, we can make an induction for
$d-l$. Observe that by the false centre theorem, by our hypothesis,
also the section of $K'$ by  
each affine $(l+1)$-plane supporting $K$ is centrally symmetric.)

In our proof of Theorem 1 (cf. \S 3), the considered sections
are in first approximation (that is, approximating
${\text{bd}}\,K^{\varepsilon }$ in the
second order, up to terms of higher order) centrally symmetric with respect to
the point of tangency (unique point of $K \cap Y$).
Therefore, for a local variant of this question
we would need to consider third order approximations of
${\text{bd}}\,K^{\varepsilon }$, of course
assuming $C^3_+$. We could suppose that for these intersections
there are some inner points (one can take, e.g., the barycentres)
such that in opposite directions the radial
functions associated to these points have values equal up to a factor
$1+O(\varepsilon ^2)$. The conclusion would be that $K^{\varepsilon }$ is
$O(\varepsilon ^2)$-close to some ellipsoid (and then $K$ would be exactly
an ellipsoid). Is this true? 

Returning to sections by all hyperplanes passing through some fixed point, 
there arises a
related question: is there a stability variant of the false centre theorem?


\head \S3 Proofs \endhead

\demo{Proof of Theorem 1 and of Remark 3}
{\bf{1.}}
It suffices to prove the case $l=d-1$. In fact, to calculate $(\partial \varrho
^t / \partial t)(u)|_{t=0}$ for $u \in S^{d-1}$, it suffices to consider some
linear $(l+1)$-subspace in ${\Cal{L}}_{l+1}$
containing $u$. There we can already calculate this quantity, using the case
of ${\Bbb{R}}^{l+1}$ and affine $l$-planes in it.

{\bf{2.}}
Hence from now on we suppose $l=d-1$.

We give the proof for Theorem 1, which concerns the $k$th intrinsic volume
$V_k( \cdot )$. However, we will always stress
(with italics or in brackets)
what properties of $V_k( \cdot )$ are used, in order to see that the proof
works also more generally for the functions $F^Y( \cdot )$ from Remark 3.

In the whole proof, when the sign $o( \cdot )$ is
applied, it is meant for $\varepsilon \to 0$.
We suppose in the whole proof that 
$\varepsilon >0$ is sufficiently small.

\newpage

Let $x \in {\text{bd}}\,K$, and let us choose a rectangular
coordinate system in ${\Bbb{R}}^d$ such that $x$
becomes the origin (thus we will have the radial function of $K$
with respect to some point of
${\text{int}}\,K$) and the lengths of the vectors in ${\Bbb{R}}^d$ in the
original polar coordinate system and this rectangular coordinate system are
identical.
Further, the hyperplane $x_d=0$ should be a tangent
hyperplane of $K$, with $K$ lying above this hyperplane. 

Then ${\text{bd}}\,K$
can be given locally, close to $x=0$, as 
$$
\cases
x_d=f(x_1, \ldots ,x_{d-1}) = f(0, \ldots , 0) + \sum _{i=1}^{d-1} f_{x_i}(0,
\ldots , 0) x_i \\
+ (1/2) \sum _{i,j=1}^{d-1} f_{x_i x_j} (0, \ldots , 0) x_i x_j +
o(x_1^2+ \ldots +x_{d-1}^2) \\
= (1/2) \sum _{i,j=1}^{d-1} f_{x_i x_j} (0, \ldots
, 0) x_i x_j + o(x_1^2+ \ldots +x_{d-1}^2)\,.
\endcases
\tag 1.1
$$
Consider the 
outer normal of $K$ at $x=0$, and the point $x^{\varepsilon} \in
{\text{bd}}\,K^{\varepsilon }$ lying on this outer normal,
below $x$ (i.e., on the negative $x_d$-axis, with $0$ included). 
Let
$$
c(x) := - (\partial / \partial t) f ^t (0, \ldots , 0) |_{t=0} \,\,(\ge 0) \,.
\tag 1.2
$$
This number $c(x)$ can be expressed by the values of the function
$[0,1] \times S^{d-1} \ni (t,u) \mapsto \varrho ^t (u) \in (0, \infty )$
and its first partial derivatives at $(0, x/ \| x \| )$. We have to use 
the transition map between two
coordinate systems, one being the polar coordinate system in ${\Bbb{R}}^d$, 
and the other one being the 
rectangular coordinate system introduced above.
However, the explicit formula is not needed. Then non-negativity
of $\partial \varrho ^t (u) / \partial t$
for $t=0$, for each $u \in S^{d-1}$, 
is equivalent to non-negativity of $c(x)$, for any boundary point
$x$ of $K$. 
From now on we will frequently write $\varepsilon $ rather than $t$, to
emphasize its smallness.
Observe that all boundary points $x^{\varepsilon }$
of $K^{\varepsilon }$ (which contains
$K$) can occur here. Namely, for $x^{\varepsilon } \in
{\text{bd}}\,K^{\varepsilon }$ we can consider its image $x$
by the nearest point map ${\Bbb{R}}^d \setminus {\text{int}}\,K 
\to {\text{bd}}\,K$: then 
$x^{\varepsilon }$ lies on the outward normal of $K$ at $x$.

Close to $x^{\varepsilon }$, ${\text{bd}}\,K^{\varepsilon }$
can be given locally as  
$$
\cases
x_d=f^{\varepsilon }(x_1, \ldots , x_{d-1})=
f ^{\varepsilon } (0, \ldots , 0) + \sum _{i=1}^{d-1} f^{\varepsilon } _{x_i}(0,
\ldots , 0) x_i \\
+ (1/2) \sum _{i,j=1}^{d-1} f^{\varepsilon } _{x_i x_j} 
(0, \ldots , 0) x_i x_j+ o(x_1^2+ \ldots +x_{d-1}^2) \\
=  -\left( c(x)+o(1) \right) \varepsilon 
+O(\varepsilon ) \sqrt{x_1^2+ \ldots +x_{d-1}^2} + \\
(1/2)
\sum _{i,j=1}^{d-1} f_{x_i x_j}(0, \ldots ,0)x_ix_j + O(\varepsilon ) (x_1^2
+ \ldots +x_{d-1}^2) \\
+ o(x_1^2+ \ldots x_{d-1}^2)\,.
\endcases
\tag 1.3
$$
In \thetag{1.3} we used that the values of the functions $f$ and
$f^{\varepsilon }$ and their
first and second derivatives with respect to $x_1, \ldots , x_{d-1}$
differ by at most $O(\varepsilon )$ (cf. \thetag{1.1}), and we used also
\thetag{1.2}.

\newpage

We will need these expansions only for the case when 
$$
x_1^2 + \ldots + x_{d-1}^2 = O(\varepsilon )
\tag 1.4
$$
(the reason for this will be given later).
So we suppose validity of \thetag{1.4} in the following.

By \thetag{1.4}, \thetag{1.3} becomes
$$
x_d = f^{\varepsilon } (x_1, \ldots , x_{d-1}) = 
-\left( c(x)+o(1) \right) 
\varepsilon + (1/2)\sum _{i,j = 1}^{d-1} f_{x_i
x_j} (0, \ldots , 0) x_ix_j \,.
\tag 1.5
$$

We have to consider the intersection $S^{\varepsilon }$
of $K^{\varepsilon }$ with the
hyperplane $x_d=0$, and have to estimate its 
$k$'th intrinsic volume, for $1 \le k \le d-1$ ($S$ for ``section''). 
For some $\varepsilon _0>0$, for all $\varepsilon \in [0, \varepsilon _0]$ we
have that $K^{\varepsilon }$ is $C^2_+$; we may suppose $\varepsilon _0 = 1$.
Therefore either $x \in {\text{int}}\,K'$, and then $S^{\varepsilon } \cap
({\text{int}}\,K') \ne \emptyset $, or $x \in {\text{bd}}\,K'$, and then
$S^{\varepsilon }$ is the one point set consisting of the point of tangency of
$K^{\varepsilon }$ with its tangent hyperplane given by $x_d=0$.
Observe that in both these cases $S^{\varepsilon }$
is the convex hull of the intersection of
${\text{bd}}\,K^{\varepsilon }$ with the hyperplane $x_d=0$.   
For $x_d = 0$, \thetag{1.5} becomes
$$
(1/2)\sum _{i,j=1}^{d-1}f_{x_i x_j} (0, \ldots ,0) x_i x_j  = 
\left( c(x)+o(1) \right) \varepsilon \,.
\tag 1.6
$$
Thus we have to estimate the $k$'th intrinsic 
volume of $S^{\varepsilon }$, which set is the convex hull of the set given
by \thetag{1.6}, that is also the star-shaped hull of the set given
by \thetag{1.6}. This star-shaped hull is
$$
S^{\varepsilon } = \{ (x_1, \ldots , x_{d-1} ) \mid
(1/2) \sum _{i,j=1}^{d-1} f_{x_i x_j}  (0, \ldots , 0)
x_i x_j \le \left( c(x)+o(1) \right) \varepsilon \} \,.
\tag 1.7
$$ 

$$
\cases
{\text{If }} c(x)>0, {\text{ then let }}
0<c_1<c(x)<c_2 {\text{ be arbitrary.}}\\
{\text{If }} c(x)=0, {\text{ then let }}
0 < c_2 {\text{ be arbitrary.}}
\endcases
\tag 1.8
$$
First we deal with the case $0 < c(x)$.
Then, for $\varepsilon >0$ sufficiently small,
we have by \thetag{1.7}
$$
\cases
S^{\varepsilon }_1 := 
\{ (x_1, \ldots , x_{d-1}) \mid \sum _{i,j=1}^{d-1}  (1/2)
f_{x_i x_j} (0,\ldots ,0) x_i x_j \le c_1 \varepsilon \} \\
\subset 
S^{\varepsilon } = 
\{ (x_1, \ldots , x_{d-1}) \mid \\
(1/2) \sum _{i,j=1}^{d-1} f_{x_i x_j} (0,\ldots ,0) x_i x_j \le 
\left( c(x)+o(1) \right) \varepsilon \} \\
\subset
S^{\varepsilon }_2 := 
\{ (x_1, \ldots , x_{d-1}) \mid (1/2) \sum _{i,j=1}^{d-1} f_{x_i x_j} (0,\ldots
,0) x_i x_j \le c_2 \varepsilon \} \,.
\endcases
\tag 1.9
$$
We are going to give lower and upper estimates for $V_k(S^{\varepsilon })$,
cf. \thetag{1.12}. For $0=c(x)$ 

\newpage

the same considerations yield only the upper
estimate for $V_k(S^{\varepsilon })$ in \thetag{1.12}:
then we use the trivial lower
estimate $0 \le V_k(S^{\varepsilon })$ rather than the one in \thetag{1.12}
(cf. \thetag{1.13}).

We write 
$$
E := \{ (x_1, \ldots , x_{d-1}) \mid (1/2) \sum _{i,j=1}^{d-1} f_{x_i x_j}
(0,\ldots ,0) x_i x_j \le 1 \} \,.
\tag 1.10
$$
Then $E$ is the convex hull of the Dupin indicatrix of $K$ at 
$x \in {\text{bd}}\,K$, taken in our chosen rectangular coordinate system.
By positive definiteness of the quadratic form $(1/2) \sum _{i,j=1}^{d-1}
f_{x_i x_j} (0,\ldots ,0) x_i x_j$, the set $E$ is an $0$-symmetric ellipsoid.
Then \thetag{1.9} can be rewritten as
$$
(c_1 \varepsilon )^{1/2} E =
S^{\varepsilon }_1 \subset 
S^{\varepsilon } \subset 
S^{\varepsilon }_2 =
(c_2 \varepsilon )^{1/2} E \,.
\tag 1.11
$$
Therefore even the third, i.e., largest 
set in \thetag{1.9} (and \thetag{1.11})
has a distance at most $O(\sqrt{\varepsilon })$ from $0$, 
so that we need to consider only such points $(x_1, \ldots ,x_{d-1})$, for
which $x_1^2+ \ldots +x_{d-1}^2 = O(\varepsilon )$. This justifies the
supposition of the validity of \thetag{1.4}.

By {\it{$k$'th degree positive homogeneity and 
monotonicity of}} $V_k( \cdot )$, \thetag{1.11} implies
$$
\cases
0 \le 
(c_1 \varepsilon )^{k/2} V_k(E) = V_k\left( (c_1 \varepsilon )^{1/2} E \right)
\le V_k(S^{\varepsilon }) \le \\
V_k\left( (c_2 \varepsilon )^{1/2} E \right) 
= (c_2 \varepsilon )^{k/2} V_k(E) \,.
\endcases
\tag 1.12
$$
As mentioned just below \thetag{1.9},
for $c(x) = 0$ we have, rather than \thetag{1.12},
$$
0 \le V_k(S^{\varepsilon }) \le (c_2 \varepsilon )^{k/2} V_k(E) \,.
\tag 1.13
$$

Hence for $0 < c(x)$ by \thetag{1.12}, while for $0 = c(x)$ by \thetag{1.13}, 
we have 
$$
V_k(S^{\varepsilon }) =
\left( \left( c(x)+o(1) \right) \varepsilon \right) ^{k/2} V_k(E) 
\tag 1.14
$$
for $\varepsilon \to 0$. 
Namely, for $0 < c(x)$ we may choose both $c_1,c_2$ arbitrarily close to $c(x)$,
and for $0 = c(x)$ we may choose $c_2 > 0$ arbitrarily close to $0 = c(x)$
(cf. \thetag{1.8}). 
We rewrite \thetag{1.14} as
$$
\lim _{\varepsilon \to 0}
V_k(S^{\varepsilon }) / [ \varepsilon ^{k/2} V_k(E) ] = c(x)^{k/2}\,.
\tag 1.15
$$
By \thetag{1.15} $K$ and 
$\lim _{\varepsilon \to 0} \left( V_k (S^{\varepsilon })/\varepsilon ^{k/2}
\right)$, for each $x \in {\text{bd}}\,K$,
determine $c(x)$ uniquely. (Recall that $k \ge 1$.)

Last, taking into account \thetag{1.2},
knowledge of this non-negative number $c(x)$, for each boundary point $x$ of
$K$, determines the non-negative partial derivative
of $\varrho ^t (u)$ 

\newpage

with respect to $t$, for $t=0$ and each $u \in S^{d-1}$.
For this we have to use
the values of the function $[0,1] \times S^{d-1} \ni (t,u) \mapsto
\varrho ^t(u) \in (0, \infty )$ and its first partial derivatives at $(0,
x/ \| x \| )$, and
use the transition map between two coordinate systems: one is the polar
coordinate system in ${\Bbb{R}}^d$, and the other one
is the rectangular coordinate system used above in the proof.
$ \blacksquare $
\enddemo


\demo{Proof of Corollary 2}
It follows immediately from Theorem 1. We only note that invariance of the first
partial derivative of the perturbation with respect to 
$t$, for $t=0$, under the map $T$
in the rectangular coordinate systems at $x, Tx \in {\text{bd}}\,K$
used in the proof of Theorem 1, 
implies its invariance in the original
polar coordinate system in ${\Bbb{R}}^d$ under the map $T$.
$ \blacksquare $
\enddemo


\demo{Proof of Theorem 4 and Remark 6}
We give the proof for Theorem 4, which concerns the $k$th intrinsic volume
$V_k( \cdot )$. However, again we will 
always stress (with italics or in brackets)
what properties of $V_k( \cdot )$ are used, in order to see that the proof
works also more generally for the functions $G^H( \cdot )$ from Remark 6.

We use the notations of the proof of Theorem 1.
In particular, in the whole proof, when the sign $o( \cdot )$ is
applied, it is meant for $\varepsilon \to 0$.
Again we suppose in the whole proof that 
$\varepsilon >0$ is sufficiently small.

{\bf{1.}}
First we consider the case $k=d$, i.e., we consider the $d$-volume of the
``caps'' cut off from $K^{\varepsilon }$ by the tangent hyperplanes of $K$.
We will write $V( \cdot )$ rather than $V_d( \cdot )$.

Till \thetag{1.8} we just use the considerations from the proof of Theorem 1.

Again, first we deal with the case $0 < c(x)$.

In \thetag{1.9} we had inclusions of (``in general'')
$(d-1)$-dimensional compact convex sets in the
hyperplane $x_d=0$. This has to be replaced by inclusions of (``in general'')
$d$-dimensional compact convex sets in ${\Bbb{R}}^d$. 

We have to investigate the ``cap'' cut off from $K^{\varepsilon }$ by the
tangent hyperplane $x_d=0$ of $K$, i.e., the set
$$
\cases
C^{\varepsilon }
:= \{ (x_1, \ldots ,x_d) \mid f^{\varepsilon }_d (x_1, \ldots ,x_{d-1}) = \\
f^{\varepsilon } (0, \ldots , 0) + \sum _{i=1}^{d-1} f^{\varepsilon
}_{x_i}(0, \ldots , 0)x_i + \\
(1/2) \sum _{i,j=1}^{d-1}f^{\varepsilon }_{x_i x_j}(0, \ldots , 0) 
x_i x_j +o(x_1^2+ \ldots + x_{d-1}^2) \le x_d \le 0\,,
\endcases
\tag 4.1
$$
($C$ for ``cap''). We are going to give lower and upper estimates for
$V(C^{\varepsilon })$, cf. \thetag{4.7}. 
For $0 = c(x)$ the same considerations
yield only \thetag{4.7} with omission of the middle expression, and
replacing the
equality sign by the $\le $ sign. Then we use the trivial lower estimate $0
\le V(C^{\varepsilon })$ rather than the one derived from \thetag{4.5}, second
inequality. However, this gives the same formula \thetag{4.7}, but omitting the
second expression there, also for $0 = c(x)$.

For $c(x)>0$ and $\varepsilon >0$ sufficiently small, 
taking into account \thetag{1.5},  
the set $C^{\varepsilon }$ in \thetag{4.1} contains

\newpage

$$
\cases
C^{\varepsilon }_1 := \{ (x_1, \ldots ,x_d) \mid -c_1 \varepsilon + \\
(1/2) \sum _{i,j=1}^{d-1} f_{x_i x_j} (0, \ldots ,0) x_i x_j 
\le x_d \le 0 \} \,,
\endcases
\tag 4.2
$$
and is contained in
$$
\cases
C^{\varepsilon }_2 := \{ (x_1, \ldots ,x_d) \mid -c_2 \varepsilon + \\
(1/2) \sum _{i,j=1}^{d-1} f_{x_i
x_j} (0, \ldots ,0) x_i x_j \le x_d \le 0 \} \,.
\endcases
\tag 4.3
$$
From \thetag{4.2} and \thetag{4.3} we have
$$
C^{\varepsilon }_1 \subset C^{\varepsilon } \subset C^{\varepsilon }_2 \,.
\tag 4.4
$$
Hence, by {\it{monotonicity and non-negativity of\,}} $V( \cdot )$, we have
$$
0 \le V(C^{\varepsilon }_1) \le V(C^{\varepsilon }) 
\le V(C^{\varepsilon }_2) \,.
\tag 4.5
$$
$$
\cases
{\text{For }} i=1,2 {\text{ the sets }} C^{\varepsilon }_i 
{\text{ are bounded by the }} (d-1){\text{-ellipsoids}} \\
S^{\varepsilon }_i {\text{ (cf. \thetag{1.9})
lying in the hyperplane given by }} x_d=0, \\
{\text{and by portions of the elliptic paraboloids given by }} x_d =  \\
-c_i \varepsilon + 
(1/2) \sum _{i,j=1}^{d-1} f_{x_i x_j} (0, \ldots ,0) x_i x_j \,,
{\text{ lying below }} S^{\varepsilon }_i \,.
\endcases
\tag 4.6
$$
The volume of $C^{\varepsilon }_i$ can be calculated as
$2/(d+1)$ times the volume of its {\it{circumscribed right cylinder 
$(C^{\varepsilon }_i)'$, with upper base the ellipsoid $S^{\varepsilon }_i$ 
and height $c_i \varepsilon $.}} 
In fact, by applying an affinity we may suppose
$E=B^{d-1}$ and $c_i \varepsilon =1$, and then we calculate
$V(C^{\varepsilon }_i)/V((C^{\varepsilon }_i)')$ by using polar
coordinates in the hyperplane given by $x_d=0$. (For the case of
$G^H( \cdot )$ we observe that
this affinity can be factorized as
the product of a magnification in a positive ratio, and a volume-preserving
affinity.)
Then for $\varepsilon >0$ sufficiently small, analogously to \thetag{1.14}
we have, both for $0<c(x)$ and for $0=c(x)$ (for $0=c(x)$ omitting the middle
term)
$$
\cases
V(C^{\varepsilon })= \left[ \kappa _{d-1}/ 
[ {\text{det}} \left( (1/2) \left( f_{x_i x_j}(0, \ldots ,0)
\right. \right. \right. \\
\left. \left. \left.
/\left( \left( c(x)+o(1) \right) 
\varepsilon \right) \right) \right) ] ^{1/2} \right] 
\cdot \left( \left( c(x)+o(1) \right) \varepsilon \right) \cdot 2/(d+1) = \\
\left[ \kappa _{d-1}/
[ {\text{det}} \left( (1/2) f_{x_i x_j}(0,0, \ldots ,0) \right) ]
^{1/2} \right] 
\cdot [2/(d+1)] \cdot \\
\left( \left( c(x)+o(1) \right) \varepsilon \right) ^{(d+1)/2} \,.
\endcases
\tag 4.7
$$
(For Remark 6 we do not have \thetag{4.7}, since we do not know
$G^H(E)$. First suppose 

\newpage

$0 < c(x)$. Then
we have the analogue of \thetag{4.5} with $G^H$, and 
by \thetag{4.6} and \thetag{1.8} we get $G^H\left( (C_1^{\epsilon })' \right) /
G^H\left( (C_2^{\epsilon })' \right) = 1 + o(1)$. Therefore, we have the
analogue of \thetag{4.7}, namely $G^H ( C^{\epsilon } ) = G^H (E') \cdot
\left( \left( c(x)+ o(1) \right) \varepsilon \right) ^{(d+1)/2}$.
Here $E'$ is the set bounded by $E$ (defined in \thetag{1.10})
and by the portion of the elliptic paraboloid given by $x_d = -1 + 
(1/2) \cdot \sum _{i,j=1}^{d-1} f_{x_i x_j}(0, \ldots , 0) x_i x_j$,
lying below $E$. The case $0 = c(x)$ is 
treated like in the case of $V( \cdot
)$, only using $G^H( \cdot )$ rather than $V( \cdot )$.)

Hence, analogously to \thetag{1.15}, we have from \thetag{4.7}
$$
\cases
\lim _{\varepsilon  \to 0} V(C^{\varepsilon }) \cdot
[ {\text{det}} \left( (1/2) \left( f_{x_i x_j}(0, \ldots
,0) \right) \right) ] ^{1/2} / \\
\left[ \kappa _{d-1} [2/(d+1) ] \varepsilon ^{(d+1)/2} \right] = c(x) ^{(d+1)/2}
\,.
\endcases
\tag 4.8
$$
By \thetag{4.8} $K$ and 
$\lim _{\varepsilon \to 0} \left(
V (C^{\varepsilon })/\varepsilon ^{(d+1)/2} \right) $, 
for each $x \in {\text{bd}}\,K$, determine $c(x)$ uniquely.
(For Remark 6, we use here once more $G^H(E')$, like below \thetag{4.7}.)
Then the 
repetition of the last paragraph of the proof of Theorem 1 finishes the
proof for the case $k=d$.

{\bf{2.}}
We turn to the case $1 \le k \le d-1$. Again, first we deal with the case $0 <
c(x)$.

For $0 = c(x)$ 
the same considerations yield for $V_k(C^{\varepsilon })$ only \thetag{4.17}.
Then we use the trivial lower estimate \thetag{4.11} rather than \thetag{4.10}.
However, this gives the same formula \thetag{4.18}, also for $0 = c(x)$.

By \thetag{1.11} 
we have $(c_1 \varepsilon )^{1/2} E = S^{\varepsilon }_1 
\subset S^{\varepsilon } \subset C^{\varepsilon }$. 
Here for $\varepsilon \to 0$ we may choose $0 < c_1 = c(x)+o(1)$
(cf. \thetag{1.8}), hence
$C^{\varepsilon } \supset \left( \left( c(x) + o(1) \right)) 
\varepsilon \right) ^{1/2} E $.
Thus, by {\it{$k$'th degree positive homogeneity and monotonicity of}}
$V_k( \cdot )$, we have  
$$
V_k(C^{\varepsilon }) / [V_k(E) \cdot \varepsilon ^{k/2}] \ge \left(
c(x)+o(1) \right) ^{k/2} \,.
\tag 4.10
$$

For $0 = c(x)$ we use the trivial estimate 
$$
V_k(C^{\varepsilon }) / [V_k(E) \cdot \varepsilon ^{k/2}] \ge 0 \,.
\tag 4.11
$$

We turn to the upper estimate. Recall that $c_2 > 0$ (cf. \thetag{1.8}).
With $(C^{\varepsilon }_2)'$ as defined below \thetag{4.6} we have by 
\thetag{4.4}
$C^{\varepsilon } \subset C^{\varepsilon }_2 \subset (C^{\varepsilon }_2)'$.
By {\it{monotonicity of}} $V_k( \cdot )$ this implies
$V_k (C^{\varepsilon }) / V_k \left( (C^{\varepsilon }_2)' \right) \le 1 $.
Equivalently, using {\it{$k$'th degree positive homogeneity of}}
$V_k( \cdot )$, we have
$$
V_k \left( C^{\varepsilon } / (c_2 \varepsilon )^{1/2} \right) /
V_k \left( (C^{\varepsilon }_2)' / (c_2 \varepsilon )^{1/2} \right) \le 1 \,.
\tag 4.14
$$
Now observe that $( C^{\varepsilon }_2)' / (c_2
\varepsilon )^{1/2} $ is a right cylinder with base $S^{\varepsilon
}_2/ (c_2 \varepsilon )^{1/2}= E$ (cf. \thetag{1.11})
and height $(c_2 \varepsilon )
/ (c_2 \varepsilon )^{1/2} = (c_2 \varepsilon )^{1/2}$. Therefore 
$( C^{\varepsilon }_2)' / (c_2 \varepsilon )^{1/2} $ 
is in the $(c_2 \varepsilon )^{1/2}$-neigh\-bour\-hood of $E$. 

Therefore we have, by {\it{continuity of\,}} 
$V_k( \cdot )$, that

\newpage

$$
V_k \left( C^{\varepsilon }_2)' / (c_2 \varepsilon )^{1/2}
\right) / V_k(E) \le 1+o(1) \,.
\tag 4.15
$$ 
Multiplying \thetag{4.14} and \thetag{4.15} we get
$$
V_k \left( C^{\varepsilon } / 
(c_2 \varepsilon )^{1/2} 
\right) / V_k(E) \le 1+o(1) \,.
\tag 4.16
$$
Since we may choose $c_2 >0$ arbitrarily close to $c(x)$ (cf. \thetag{1.8}), 
therefore once more by {\it{$k$'th degree positive homogeneity of}}
$V_k( \cdot )$, \thetag{4.16} implies
$$
V_k(C^{\varepsilon }) / [V_k(E) \cdot \varepsilon ^{k/2}] \le c_2^{k/2} \left(
1+o(1) \right) = \left( c(x)+o(1) \right) ^{k/2} \,.
\tag 4.17
$$

Then for $0 < c(x)$ \thetag{4.10}, while 
for $0=c(x)$ \thetag{4.11}, together with
\thetag{4.17} give
$$
V_k(C^{\varepsilon }) / [V_k(E) \cdot \varepsilon ^{k/2}] = 
\left( c(x)+o(1) \right) ^{k/2} \,.
\tag 4.18
$$
We rewrite \thetag{4.18} as
$$
\lim _{\varepsilon \to 0} 
V_k(C^{\varepsilon }) / [V_k(E) \cdot \varepsilon ^{k/2}] = c(x)^{k/2} \,.
\tag 4.19
$$
By \thetag{4.19} $K$ and 
$\lim _{\varepsilon \to 0} \left(
V_k(C^{\varepsilon })/\varepsilon ^{k/2} \right) $, 
for each $x \in {\text{bd}}\,K$,
determine $c(x)$ uniquely. (Recall that $k \ge 1$.)

Then the repetition of the last paragraph of the proof of Theorem 1 finishes the
proof for the case $1 \le k \le d-1$.
$ \blacksquare $
\enddemo


\demo{Proof of Corollary 5}
It follows immediately from Theorem 4, also taking into account the proof of
Corollary 2.
$ \blacksquare $
\enddemo



\Refs

\widestnumber\key{WWW}



\ref
\key APR
\by P. W. Aitchison, C. M. Petty, C. A. Rogers
\paper A convex body with a false centre is an ellipsoid
\jour Mathematika
\vol 18
\yr 1971
\pages 50-59
\MR {\bf{45\#}}{\rm{2572}}. 
\endref 

\ref
\key BL
\by J. A. Barker, D. G. Larman
\paper Determination of convex bodies by certain sets of sectional volumes
\jour Discrete Math.
\vol 241
\yr 2001
\pages 79-96
\MR {\bf{2002j:}}{\rm{52001}}. 
\endref 

\ref
\key F 
\by K. Falconer
\paper $X$-ray problems for point sources
\jour Proc. London Math. Soc. (3)
\vol 46
\yr 1983
\pages 241-262
\MR {\bf{85g:}}{\rm{52001a}}. 
\endref 

\ref
\key Ga83 
\by R. J. Gardner
\paper Symmetrals and $X$-rays of planar convex bodies
\jour Arch. Math. (Basel)
\vol 41
\yr 1983
\pages 183-189
\MR {\bf{85c:}}{\rm{52007}}. 
\endref 

\ref
\key  Ga 
\book Geometric Tomography; Second edition
\by R. J. Gardner
\publ Encyclopedia of Math. and its Appl. {\bf{58}}, Cambridge Univ. Press
\publaddr Cambridge 
\yr 1995; 2006
\MR {\bf{96j:}}{\rm{52006}}; {\bf{2007i:}}{\rm{52010}}.
\endref 

\ref
\key GK
\by V. Gorkavyy, D. Kalinin
\paper Barker-Larman problem in the hyperbolic plane
\jour Results Math.
\vol 68
\yr 2015
\pages 519-525
\MR {\bf{3407571}}{\rm{}}. 
\endref 

\newpage

\ref
\key Gr 
\book Geometric Applications of Fourier Series and Spherical Harmonics
\by H. Groemer
\publ Encyclopedia of Math. and its Appl. {\bf{61}}, Cambridge Univ. Press
\publaddr Cambridge
\yr 1996
\MR {\bf{97j:}}{\rm{52001}}.
\endref 

\ref
\key H
\by P. C. Hammer
\paper Diameters of convex bodies
\jour Proc. Amer. Math. Soc.
\vol 5
\yr 1954
\pages 304-306
\MR {\bf{15,}}{\rm{819b}}. 
\endref 

\ref
\key J-CM
\by J. Jer\'onimo-Castro, T. B. McAllister
\paper Two characterizations of ellipsoidal cones
\jour J. Convex Anal.
\vol 20
\yr 2013
\pages 1181-1187
\MR {\bf{3184302}}{\rm{}}. 
\endref

\ref
\key L
\by D. G. Larman
\paper A note on the false centre theorem 
\jour Mathematika
\vol 21
\yr 1974
\pages 216-227
\MR {\bf{50\#}}{\rm{14490}}. 
\endref 

\ref
\key LM-A
\by D. G. Larman, E. Morales-Amaya
\paper On the false pole problem
\jour Monatsh. Math.
\vol 151
\yr 2007
\pages 271-286
\MR {\bf{2008d:}}{\rm{52006}}. 
\endref 

\ref
\key MM 
\by E. Makai, Jr., H. Martini
\paper Centrally symmetric convex bodies and sections having maximal
quermassintegrals
\jour Studia Sci. Math. Hungar.
\vol 49
\yr 2012
\pages 189-199
\MR {\bf{3058386}}.
\endref 

\ref
\key MM\'O
\by E. Makai, Jr., H. Martini, T. \'Odor
\paper Maximal sections and centrally symmetric convex bodies 
\jour Mathematika
\vol 47
\yr 2000
\pages 17-30
\MR {\bf{2003e:}}{\rm{52005}}.
\endref 

\ref
\key MM-A
\by L. Montejano, E. Morales-Amaya
\paper Variations of classic characterizations of ellipsoids and a short proof
of the false centre theorem
\jour Mathematika
\vol 54 {\rm{(1-2)}}
\yr 2007 
\pages (1-2), 35-40
\MR {\bf{2009c:}}{\rm{52008}}. 
\endref 

\ref 
\key RYZ
\by D. Ryabogin, V. Yaskin, A. Zvavitch
\paper Harmonic analysis and uniqueness questions in convex geometry
\jour 
In: Recent Advances in Harmonic Analysis and Applications (Eds. D. Bilyk et
al.), Springer, New York, 2013, Proc. Math. Stat. {\bf{25}} 
\pages 327-337  
\MR {\bf{3066896}}.
\endref   

\ref
\key Sa 
\by L. A. Santal\'o
\paper Two characteristic properties of circles on a spherical surface
(Spanish)
\jour Math. Notae
\vol 11
\yr 1951
\pages 73-78
\MR {\bf{14,}}{\rm{495c}}. 
\endref 

\ref
\key Sch 
\book Convex Bodies: the Brunn-Minkowski Theory; Second expanded edition
\by R. Schneider
\publ Encyclopedia of Math. and its Appl., {\bf{44}}; {\bf{151}},
Cambridge Univ. Press
\publaddr Cambridge
\yr 1993; 2014
\MR {\bf{94d:}}{\rm{52007}}; {\bf{3155183}}.
\endref 

\ref
\key So
\by V. Soltan
\paper Convex hypersurfaces with hyperplanar intersections of their
homothetic copies
\jour J. Convex Anal.
\vol 22 
\yr 2015
\pages 145-159
\MR {\bf{3346184}}{\rm{}}. 
\endref 

\ref
\key XMCh
\by G. Xiong, Y.-w. Ma, W.-s. Cheung
\paper Determination of convex bodies from $\Gamma $-section functions
\jour J. Shanghai Univ.
\vol 12 {\rm{(3)}}
\yr 2008
\pages 200-203
\MR {\bf{2009d:}}{\rm{52014}}. 
\endref 

\ref
\key Y
\by V. Yaskin
\paper Unique determination of convex polytopes by non-central sections
\jour Math. Ann.
\vol 349
\yr 2011
\pages 647-655
\MR {\bf{2012a:}}{\rm{52011}}. 
\endref 

\ref
\key YZh
\by V. Yaskin, N. Zhang
\paper Non-central sections of convex bodies
\jour http://arxiv.org/abs/1509.
\newline
08174
\vol
\yr
\pages 
\endref 




\endRefs
\enddocument